\documentclass{article}
\usepackage{graphicx} 
\usepackage{amsmath}
\usepackage{amsfonts,amssymb,amsthm,enumitem}
\usepackage{mathtools}
\usepackage{nicefrac,setspace}
\usepackage{xcolor}
\usepackage{tikz}
\usepackage[framemethod=tikz]{mdframed}
\usepackage{pgfplots}
\usepackage{subfig}
\usepackage{graphicx,float}
\pgfplotsset{compat=1.18}
\usepackage[font=footnotesize,labelfont=bf]{caption}
\usepackage{hyperref}       
\usepackage{url}
\usepackage{diagbox}
\newtheorem{definition}{Definition}

\begin{document}
\title{A practical, fast method for solving sum-of-squares problems for very large polynomials}
\title{A practical, fast method for solving sum-of-squares problems for very large polynomials}

\author{Daniel Keren \and Margarita Osadchy \and Roi Poranne\thanks{University of Haifa}}

\date{August 2024}
\maketitle

\begin{abstract}
Sum of squares (SOS) optimization is a powerful technique for solving problems
where the positivity of a polynomials must be enforced. The common approach to solve an SOS problem is by relaxation to a Semidefinite Program (SDP). The main advantage of this transormation is that SDP is a convex problem for which efficient solvers are readily available. However, while considerable progress has been made in recent years, the standard approaches for solving SDPs are still known to scale poorly.
Our goal is to devise an approach that can handle larger, more complex problems than is currently possible.

The challenge indeed lies in how SDPs are commonly solved. State-Of-The-Art approaches rely on the interior point method, which requires the factorization of large matrices.
We instead propose an approach inspired by polynomial neural networks, which exhibit excellent performance when optimized using techniques from the deep learning toolbox.
In a somewhat counter-intuitive manner, we replace the convex SDP formulation with a non-convex, unconstrained, and \emph{over parameterized} formulation, and solve it using a first order optimization method.
It turns out that this approach can handle very large problems, with polynomials having over four million coefficients, well beyond the range of current SDP-based approaches.
Furthermore, we highlight theoretical and practical results supporting the experimental success of our approach in avoiding spurious local minima, which makes it amenable to simple and fast solutions based on gradient descent.
In all the experiments, our approach had always converged to a correct global minimum, on general (non-sparse) polynomials, with running time only slightly higher than linear in the number of polynomial coefficients, compared to higher than quadratic in the number of coefficients for SDP-based methods. 
\end{abstract}
\section{Introduction}
Determining whether a polynomial is positive inside a domain is a problem that appears often in various application, ranging from graph theory and combinatorics \cite{9fef820b69d243f2a501e933b30bd977}, geometric modeling and robotics,\cite{Keren2024,DBLP:journals/pami/KerenG99,DBLP:journals/tog/MarschnerZPS21,DBLP:conf/siggraph/ZhangM0T23,ahmadi2017geometry},
and optimization \cite{DBLP:journals/siamjo/Lasserre01,Blekherman2012,9fef820b69d243f2a501e933b30bd977}.
Unfortunately, the problem is known to be NP-Complete in general, even for degree 4 (quartics) \cite{9fef820b69d243f2a501e933b30bd977}.
Sum Of Squares programming (SOS), is a popular relaxation, in which the positivity constraint is replaced by the requirement that the polynomial in question be a sum of squares of polynomials.
This problem in particular has inspired a decades long body of research, both in pure mathematics \cite{Reznick,Powers2021,putinar1993} and optimization \cite{DBLP:journals/siamjo/Lasserre01,parrilo2000,Blekherman2012}, and many  software packages were developed to this end \cite{sostools, seiler2013sosopt}.
Most of these rely on transforming the problem into an SDP, which can be solved fairly reliably.
However, despite considerable progress on SDP solvers in the last decades, they are notorious for their poor scalability, both in terms of time and memory requirements.
This is further exacerbated by the rate of growth in complexity of SOS problems, as the number of dimensions and degree increases.
As a result, there is somewhat of an informal upper bound on the size of the problem can be feasibly solved on a standard machine.
Our goal is to circumvent this limitation and propose an alternative approach.

The challenge lies it that, while SDP do indeed scale poorly, an SDP is a convex problem, and thus it is guaranteed to find a solution.
Other formulations are generally not convex.
However, it turns out the the formulation we propose herein is pseudoconvex, and therefore gradient based methods do reach a global minimum.
Our formulation is inspired by recent \emph{polynomial} neural networks \cite{}.
As far as we know, there is no agreed upon definition of a polynomial neural network, but it should be understood as a neural network that outputs polynomials, for example, in the case where activation function are polynomial.
More specifically, we design a network that outputs an SOS polynomial by construction, and given a target polynomial, ask whether that polynomial can be the output of our network.
A positive answer means that the target polynomial is an SOS.
In contrast to common wisdom in the realm of SOS problems, we utilize a redundant, \emph{over parameterized} model.
We demonstrate that this redundancy in fact contributes to faster convergence.
\subsection{Previous Work}
\label{sec:prev}
A widely used method for determining whether a polynomial in $n$ variables and with degree $2d$, $p(x_1 \ldots x_n)$, is an SOS, is to use the well-known fact (\cite{Powers2021}) that $p$ is a SOS iff it satisfies the following identity for some PSD matrix $B$:
\begin{equation} \label{eq:basicsos}
    p(x_1 \ldots x_n) = m_dBm_d^t
\end{equation}
where $m_d$ is a vector consisting of all monomials of total degree $ \leq d$; for example, for $n=d=2$, 
\[
m_d = (x_1^2,x_1x_2,x_2^2,x_1,x_2,1)
\]
Hereafter we will restrict attention to {\em homogeneous} polynomials (often referred to as {\em forms}), for which the total degree of every monomial is exactly $2d$ (e.g. for $n=d=2$, $m_d = (x_1^2,x_1x_2,x_2^2)$. It 
follows from homogenization \cite{Powers2021} that the problem for general polynomials in $n$ variables is naturally identical to that with homogeneous polynomials in $n+1$ variables, hence no generality is lost.

Eq. \ref{eq:basicsos} is a semidefinite programming (SDP) feasibility problem \cite{DBLP:journals/siamjo/Lasserre01,9fef820b69d243f2a501e933b30bd977}, 
with a matrix of size ${n+d-1 \choose n-1} \times {n+d-1 \choose n-1}$, 
and ${n+2d-1 \choose n-1}$ constraints (this follows from the fact
that there are ${k+d-1 \choose k-1}$ homogeneous 
monomials of degree $d$ in $k$ variables). 
While convex, and in principle approachable with interior point methods, the problem becomes quite difficult for large $n,d$; for example, existing SDP algorithms cannot handle a general quartic (fourth degree) polynomial in 40 variables ($n=40, d=2$); they either abort due to memory allocation problems, or are very slow 
\cite{DBLP:journals/corr/abs-2205-11466,DBLP:conf/nips/BoumalVB16,dsossdsos}. A summary and analysis of the SDP problem complexity as a function of $n,d$ can be found in \cite{DBLP:journals/algorithmica/JiangNW23}. Asymptotically, for quartic polynomials in $n$ variables, the running time is
$O(n^{9.5})$.

The inability of interior point methods to handle SOS problems for large $n$ had inspired the following directions of research: 
\begin{itemize}
\item {\em Sparse} polynomials, i.e. with relatively large $n,d$, but with the majority of coefficients equal to zero. Some recent work in this vein is reported in \cite{DBLP:conf/issac/WangLX19,DBLP:journals/siamjo/WangML21a}.
\item Restricting the SOS polynomials to be of a special type, for example $m_dBm_d^t$ with $B$ diagonally dominant. While allowing to tackle large $n$, these 
methods are not guaranteed to find the correct solution, as they restrict the
search space
\cite{dsossdsos}.
\item Directly solving the problem in Eq. \ref{eq:basicsos}, as follows. Using the fact that a symmetric matrix $B$ is PSD iff it can be written as $PP^t$ for some matrix $P$, one may try to minimize 
\begin{equation} \label{eq:nonconvex}
||p(x) - m_dPP^tm_d^t||^2
\end{equation}
where $x = (x_1 \ldots x_n)$, and the difference between the polynomials $p(x)$ and $m_dPP^tm_d^t$ is understood to be the squared norm of the difference between their coefficient vectors. It is this paradigm which we follow here. 
\end{itemize}
The problem in Eq. \ref{eq:nonconvex} is non-convex in $P$. Some work in this vein is presented in \cite{DBLP:journals/siamjo/PappY19}, in which the SDP approach is circumvented by directly optimizing over the cone of SOS polynomials, using non-symmetric conic optimization. In \cite{DBLP:journals/oms/BertsimasFS13}, an accelerated first-order method is applied to the SOS problem. A recently introduced approach, which achieves very fast running times for univariate polynomials with extremely high degrees, is presented in \cite{DBLP:journals/corr/abs-2205-11466}, in which is it proved that, for polynomials with a single variable, the corresponding SOS problem has no spurious local minima. The fact that every univariate SOS polynomial can be represented as a sum of only 
{\em two} squares, plays an important role in the solution.

We next elaborate on a general result concerning the problem of spurious local minima in a general type of SDP problems. 
\subsection{The Burer-Monteiro approach}
\label{sec:bm}
In a landmark series of papers \cite{DBLP:journals/mp/BurerM03,DBLP:journals/mp/BurerM05,DBLP:journals/siamjo/BurerMZ02}, Burer and Monteiro demonstrated that, in many important cases, the inherent difficulty in solving SDP problems can be substantially reduced by noting that the
problem in Eq. \ref{eq:nonconvex} has no spurious local minima; followup work, with a careful analysis of their approach, is presented in
\cite{DBLP:conf/nips/BoumalVB16,910ba8f7f19d4cddb1f7d0d0fa574105}. 
 Specifically, the Burer-Monteiro (BM hereafter) paradigm is typically applied when the number of equality constraints in 
 the SDP problem, denoted $m$, is relatively small. 
 In that case, the sought $l \times l$ PSD matrix 
can be replaced by 
$YY^t$, where $Y$ is of size $l \times p$, as long as
$p(p+1) \geq 2m$, and the resulting problem -- while non-convex -- will not suffer from spurious local minima.

This paradigm offers an attractive solution when the
number of constraints is small, allowing to replace the $l^2$ variables of the $l \times l$ PSD matrix by $2lp$ variables, where $p$ is of the order of magnitude of the square root of the number of constraints.
For the problem of non-sparse polynomials, as described in the Introduction, the number of constraints is very large, as it equals the number
of coefficients, so $m={n+2d-1 \choose n-1}$, while 
$l={n+d-1 \choose n-1}$. For example, for the largest
polynomial problem we solve in this paper ($n=100,d=4$) we have $l=5,050, m=4,421,275$. This
entails a value of $p=\sqrt{2m} \approx 2,973$, and
an overall number of variables $lp \approx 1.501 \cdot 10^7$; 
however, the standard parameterization of the PSD matrix as a
square of a symmetric $5,050 \times 5,050$ matrix, 
entails $1.275 \cdot 10^7$ parameters, which is {\em less}
than what the BM approach requires.

Still, we argue that the rarity of spurious local minima makes it worthwhile to forgo the SDP paradigm, and solve an
unconstrained problem with {\em more} variables, as this avoids the very difficult SDP-based solution. This
is supported by the experimental results, which are
presented following the description of the optimization procedure.
\subsection{The Pythagoras number of the polynomial ring, and  over-parameterization for the SOS problem}
\label{sec:over}
To recap, the BM approach seeks to replace the minimization of $\phi(X)$, where $X$ ranges over $l \times l$ PSD matrices, and $\phi$ is convex, with the minimization of $\phi(YY^t)$, where $Y$ is an unconstrained matrix of size $l \times p$, and $k$ is at least $\sqrt{2m}$, where $m$ is the number of constraints. While the size of the problem increases with $k$, it was noted that for the above problem, as well as for neural networks with polynomial activation
\cite{DBLP:journals/corr/abs-2207-01789,pmlr-v80-du18a}, increasing $k$ helps in avoiding spurious local minima. For the case of SOS polynomials, the minimal rank of $Y$ is related to the {\em Pythagoras number} of the ring of polynomials $R$ (in our case, $R$ is the set of quartic polynomials with $n$ variables, over the reals). The Pythagoras number of a ring $R$ is defined as follows:
\begin{definition}
The Pythagoras number of a ring $R$, denoted $P(R)$, is the smallest positive integer $j$ which satisfies that if $r \in R$ is a sum of squares, it is a sum of $j$ squares.
\end{definition}
For the ring of quartics in $n$ variables, an upper bound for $P(R)$ is provide in \cite{Reznick}:
\begin{equation} \label{eq:pyt}
\frac{1}{2\sqrt{3}}\left(n^2+3n+1\right) - \frac{1}{2} +  o(1) 
\end{equation}

The relevance of $P(R)$ to the minimal rank of $Y$ follows from the following consideration. As in Section \ref{sec:prev}, let $m_2$ denote the vector of all monomials of degree 2 in $n$ variables. Then, every SOS polynomial can be expressed as
\[
\sum_{i=1}^{P(R)} \langle v_i,m_2 \rangle^2 = 
m_2^t \left(\sum_{i=1}^{P(R)} v_iv_i^t \right) m_2
\]
and, evidently, $\sum_{i=1}^{P(R)} v_iv_i^t$ is of rank $P(R)$. 

However, in light of the above discussion, we chose to over-parameterize the solution in order to minimize spurious local minima, and have chosen a full rank matrix, instead of $\left(\sum_{i=1}^{P(R)} v_iv_i^t \right)$ (that is, we use a sum of ${n+1 \choose 2}$ squares). Interestingly, on the average, this over-parameterization also yields convergence in fewer iterations than when using a matrix $Y$ with rank equal to the Pythagoras number (Section \ref{sec:experiments}). 

Due to the bound in Eq. \ref{eq:pyt}, our parameterization is larger by a {\em constant multiplicative factor} from the minimal rank one; as discussed in 
\cite{DBLP:journals/corr/abs-2207-01789}, this substantially assists in avoiding spurious local minima. 
Asymptotically, 
a bound on the Pythagoras number for $n$ variables
and degree $m$, is provided by 
$\sqrt{\frac{2}{m!}}n^{-\frac{m}{2}}$, and the 
over-parameterization factor can be chosen to
fit the bounds in \cite{DBLP:journals/corr/abs-2207-01789}.

\section{The network}
\label{sec:network}
\begin{figure}[h!] 
\centering
\includegraphics[scale=0.36]{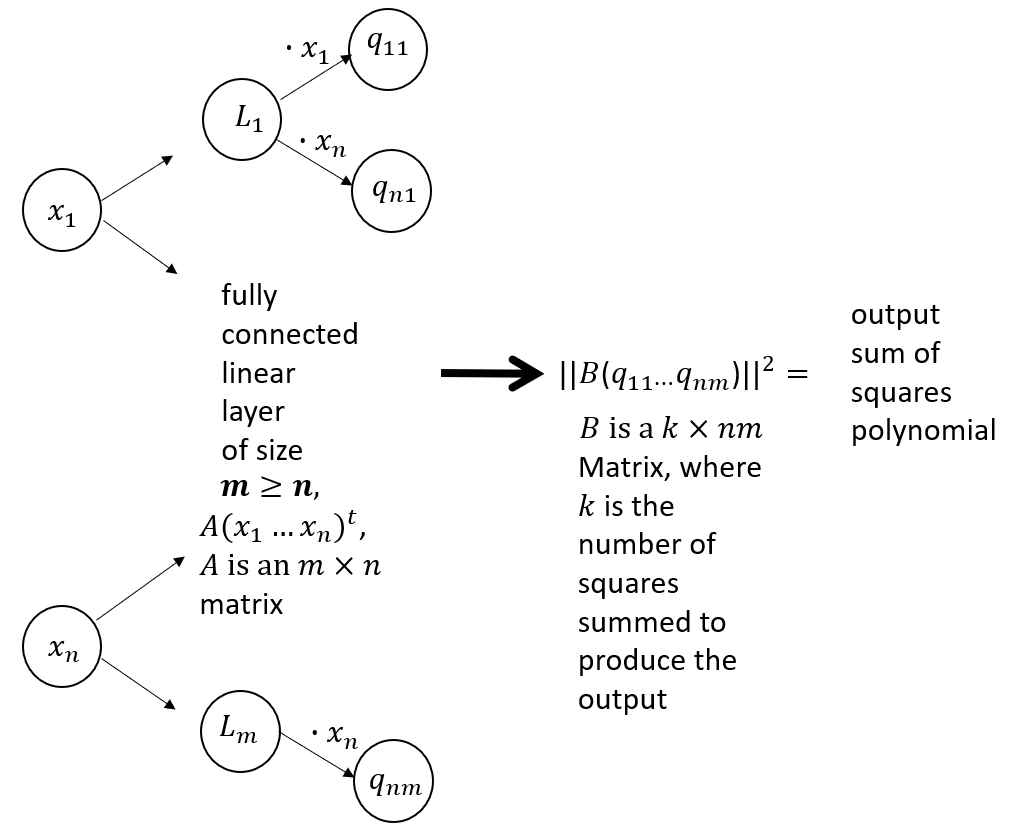}
    \caption{A schematic diagram of the
    network for generating quartic polynomials in $n$ variables. Note that it {\em symbolic}, i.e its output is not numeric, but a symbolic polynomial. The network weights are encoded in the matrices $A(m \times n),B(k \times nm)$, where $m \geq n$.}
    \label{fig:NN}
\end{figure}

We adopt a simple network for generating 
quartic polynomial in $n$ variables 
\cite{DBLP:journals/pami/ChrysosMBDPZ22,DBLP:conf/cvpr/ChrysosWDC23,DBLP:conf/nips/KileelTB19}.
Our network can be described as accepting a vector $\mathbf{x} \in \mathcal{R}^n$ and returning 
\begin{equation} \label{eq:NN}
||BZA\mathbf{x}||^2
\end{equation}
where $A,B$ are matrices (the network weights), and $Z$ is an "augmentation operation", which accepts $\mathbf{v} \triangleq A\mathbf{x}$ and outputs a vector of length $n|\mathbf{v}|$, defined by 
\((x_1\mathbf{v},x_2\mathbf{v} \ldots x_n\mathbf{v})
\). The network is schematically describe in
Fig. \ref{fig:NN}.
The network consists of a fully connected linear layer, followed by a non-linear layer, which augments the $m$ linear elements 
built by the first layer, by multiplying all
of them by $x_1 \ldots x_n$. The resulting layer (of size $mn$) is then multiplied by a
$k \times (mn)$ matrix $B$, and the output
of the network is defined as the norm squared of the result; hence, $k$ is the
number of squares summed to produce the 
output polynomial.

Note that the augmentation layer can be described by a product with a matrix composed of blocks, each of which is a product of $x_i$ with the $|\mathbf{v}| \times |\mathbf{v}|$ identity matrix. 
For example, if $n=2$, and $A$ is $3 \times 2$, then $|\mathbf{v}|=3$, and $Z$ equals
\[
Z = 
\begin{pmatrix}
x_1 & 0 & 0 \\
0 & x_1 & 0 \\
0 & 0 & x_1 \\
x_2 & 0 & 0 \\
0 & x_2 & 0 \\
0 & 0 & x_2
\end{pmatrix}
\]
In order to write down the coefficients of the 4th degree homogeneous polynomial defined by Eq. \ref{eq:NN}, we need to compute its 4th order partial derivatives for the respective monomials. For example, let us compute $\frac{\partial^4}{\partial x_1^4}$: in the following, a left-to-right arrow designates differentiation by $x_1$, and for brevity, $f$ denotes $BZA\mathbf{x}$. Note that $f$ is a (homogeneous) quadratic in $\mathbf{x}$, hence its third derivatives vanish, and that $\mathbf{x},Z$ are linear in $\mathbf{x}$, hence they have only a first non-zero derivative. To reduce equation clutter, we denote e.g. the partial derivative by $x_1$ as $f_{x_1}$.
\[
||f||^2 \rightarrow 2\langle f,f_{x_1} \rangle  \rightarrow 2(\langle f_{x_1},f_{x_1}\rangle + \langle f,f_{{x_1}{x_1}}\rangle )
\rightarrow\]
 \[2(2\langle f_{{x_1}{x_1}},f_{x_1}\rangle + \langle f_{x_1},f_{{x_1}{x_1}}\rangle ) 
\rightarrow 6\langle f_{{x_1}{x_1}},f_{{x_1}{x_1}}\rangle \]

and diving by $4!$, we see that the coefficient of $x_1^4$ equals
$\frac{1}{4}\langle f_{{x_1}{x_1}},f_{{x_1}{x_1}}\rangle$.

Proceeding in the same manner for all combination of 4 partial derivatives yields the following, where $c$ denotes the respective monomial coefficient:
\[c(x_i^4) = \frac{1}{4}\langle f_{{x_i}{x_i}},f_{{x_i}{x_i}}\rangle \]
\[c(x_i^3x_j) = \langle f_{{x_i}{x_i}},f_{{x_i}{x_j}}\rangle \]
\[c(x_i^2x_j^2) = \langle f_{{x_i}{x_j}},f_{{x_i}{x_j}}\rangle + \frac{1}{2} 
\langle f_{{x_i}{x_i}},f_{{x_j}{x_j}}\rangle \]
\[ c(x_i^2x_jx_k) = \langle f_{{x_i}{x_i}},f_{{x_j}{x_k}}\rangle +
2\langle f_{{x_i}{x_j}},f_{{x_i}{x_k}}\rangle \]
\[c(x_ix_jx_kx_l) = 2(\langle f_{{x_i}{x_j}},f_{{x_k}{x_l}}\rangle +
\langle f_{{x_i}{x_k}},f_{{x_j}{x_l}} \rangle
+ \langle f_{{x_i}{x_l}},f_{{x_j}{x_k}})\rangle \]
These identities can be understood as follows. For each set of four variables (with repetitions), by which the derivatives are to be computed, the result is obtained by summing all products of second order derivatives of $f$ by two disjoint sets whose union equals the variable set, and dividing by the product of factorials of the distinct powers in the monomial.

It remains only to compute $f_{{x_i}{x_j}}$. This equals 
\[BZ_{x_i}A\mathbf{x}_{x_j}+BZ_{x_j}A\mathbf{x}_{x_i} \]
which can be efficiently computed by noting that
$\mathbf{x}_{x_i}$ is a vector with the $i$-th coordinate 1 and all others zero, and $Z_{x_i}$ has the same structure as $Z$, but with the $i$-th block equal to the identity and the rest 0. Hence, $BZ_{x_i}A\mathbf{x}_{x_j}$ is the product of a sub-block of $B$ and a column of $A$, and its computation is much faster than computing the network's output on a single input.

These considerations directly extend to higher degrees, however, here we restrict ourselves to quartics. The problem of testing positivity for quartics is NP-complete and in addition, quartics are directly
related to important problems in graph theory and combinatorics (matrix copositivity, the partition problem, max-cut and more, \cite{9fef820b69d243f2a501e933b30bd977}).

\subsection{The first layer}
The first layer of the network we use is fully linear, and it does not enrich the set of polynomials that the network can produce. However, it proved very useful in reducing the running time in cases in which the coefficients of the input polynomial are not uniform (e.g. some are much larger than the others). More on this in Section 
\ref{sec:experiments}. 

\section{Experimental results}
\label{sec:experiments}
The proposed algorithm was tested on quartic polynomials with a number of variables ranging from 10 (715 coefficients) to 100 (4,421,275 coefficients). In all cases tested, when the input polynomial was a SOS, the error was very small, with the norm of the error vector equal to roughly $10^{-8}$ of the coefficient vector's norm. The error is defined by Eq.
\ref{eq:nonconvex}.

All tests were run on an Intel i7-6500U CPU laptop, clocked at 2.50GHz, with 16GB RAM. The code was written in Python/Jax with automatic differentiation, and optimization was carried out using the basic ADAM approach, w/o any parameter optimization. The initial guess was randomly chosen, and did not affect convergence. All results are averaged over 50 runs. While the success of the algorithm is measured
by the minimal error over all the iterations, we have
plotted the error value for the entire optimization process. The error
is un-normalized by the number of coefficients.
\subsection{Comparison with the interior point method}
\label{sec:sdp}
For comparison, the input polynomials were also tested with an SDP, interior point method, in which the polynomial is represented as $m_2Qm_2^t$, where $m_2$ is the vector of 2nd degree monomials, and $Q$ is a PSD matrix. There is no objective function, and the problem is a feasibility-only one, namely: does there
exist an SPd $Q$ yielding the coefficients of the input polynomial?

For both algorithms, we have excluded from the running time the part in which the coefficients are extracted and the optimization problem defined (as it is identical for all inputs). The inputs were generated by $m_2Qm_2^t$, where $Q$ is a PSD matrix defined by $BB^t$, and $B$'s elements are random uniform in
$[-1,1]$. For large enough $n$, the SDP approach complexity was on par with
results reported in the literature, $O(n^{9.5})$, while the complexity of our approach was roughly 
$O(n^{4.5}$. We have tested the SDP approach using the current
CVXPY package \cite{diamond2016cvxpy} only up to $n=30$, since for larger values, the "preparation"
stage, which consists of defining the problem's parameters, is too long (376,000 seconds for $n=30$). However, we note that the "net" time  ($O(n^{9.5}$) for
SDP run with 100 variables, will be inconceivably high; actually, we have not been able to find in the literature results for
applying the general SDP approach, for general, non-sparse quartics, beyond 30 variables.
~~\\
~~\\
\begin{tabular}{|c|c|c|}
  \hline
  \diagbox{Number of Variables}{Running Time (seconds)} & SDP & Our approach \\
  \hline
  10 (715 coefficients) & 2.7 & 32.7 \\
  \hline
   15 (3,060 coefficients)& 32 & 46 \\
  \hline
     20 (8,855 coefficients) & 255 & 89 \\
  \hline
     25 (20,475 coefficients)& 1,990 & 234 \\
  \hline
     30 (40,920 coefficients) & 9,266 & 542 \\
  \hline
     100 (4,421,275 coefficients)  & NA & 151,000 \\
  \hline
\end{tabular}
\subsection{Different ranks of the $B$ matrix}
As discussed in Section \ref{sec:over}, it turned out to be advantageous to over-parameterize the matrix $B$, which
defines the number of squares summed to yield the output
polynomial; this not only results in the algorithm avoiding
spurious local minima, but reduces the number of iterations
required for convergence. 

\begin{figure}[h!]
\centering
\includegraphics[scale=0.36]{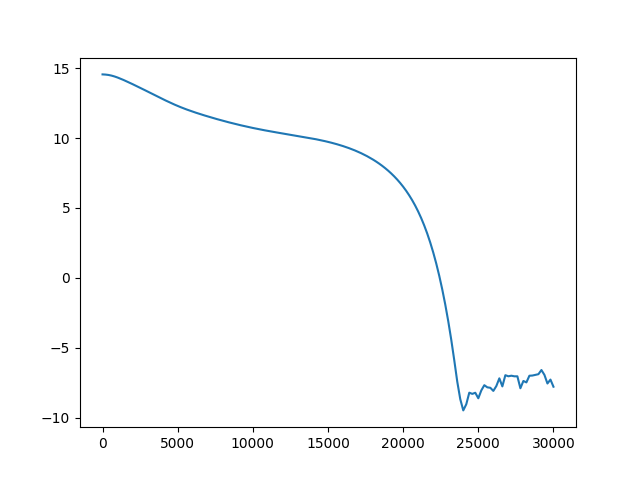}
\includegraphics[scale=0.36]{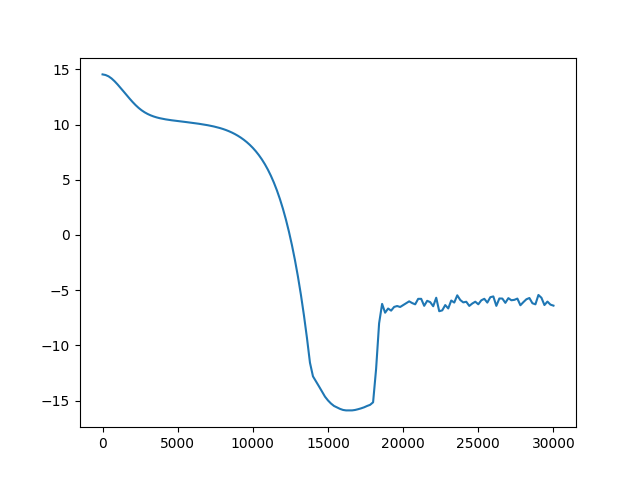}\\
\includegraphics[scale=0.36]{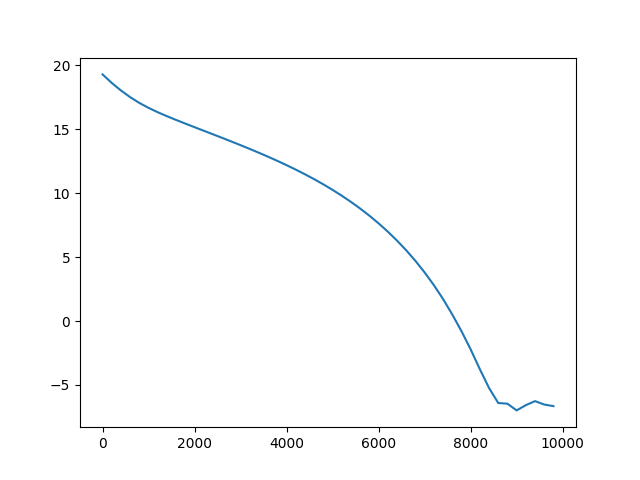}
    \caption{Top left: log of the error (vertical) vs. iteration number (horizontal), for $B$ with rank 110 (the upper bound on the Pythagoras number, Section \ref{sec:over}), vs. (top right) $B$ with rank 465
    (the number of second degree monomials with 30
    variables). Bottom: convergence for 100 variables. Note that the number of required iterations is on
    the same order of magnitude as for 30 variables.}
\end{figure}
\subsection{The effect of the first layer for polynomials with non-uniform coefficients}
Typically, SOS polynomials are sampled from a uniform distribution \cite{DBLP:journals/siamjo/WangML21a,DBLP:conf/issac/WangLX19,DBLP:journals/corr/abs-2205-11466}. However,
when running without the first layer (alternatively, setting $A$ in Section \ref{sec:network} to the $n \times n$ identity matrix), considerably slows down convergence when the
coefficients are not uniform. Below we compare results with $A=I$ vs. $A$ = a general matrix, for
"uniform coefficients" (as defined in Section \ref{sec:sdp}), but with the coefficient of $x_i^4$ set to 10,000, for some random $i$.
\begin{figure}[h!]
\centering
\includegraphics[scale=0.36]{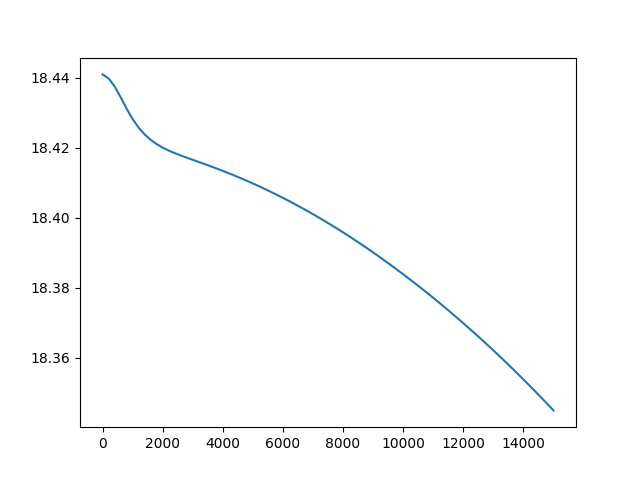}
\includegraphics[scale=0.36]{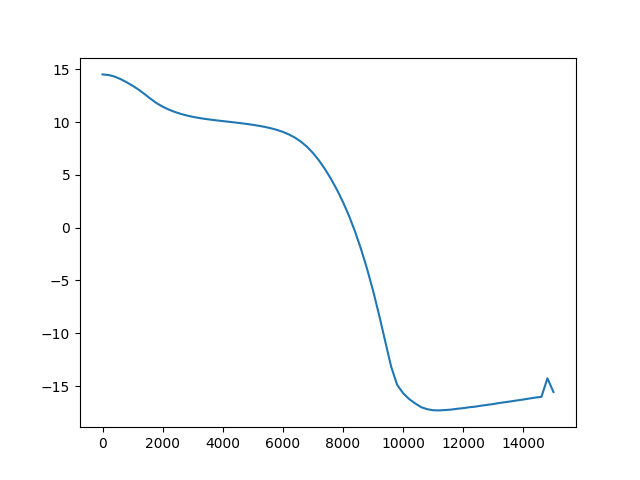}
    \caption{Left: error vs. number of iterations with $A=I$ and
    one coefficient set to $10^4$. Right: same, for general $A$.}
\end{figure}
We have also applied a more general perturbation, in which all
coefficients of $x_i^4$ and $x_i^2x_j^2$ were increased by
a random amount, uniformly distributed in $[0,10^4]$. The 
application of $A$ makes a crucial difference (although
convergence is slower than for one perturbed coefficient); and when $A$ is
also over-parameterize (taken to be not $n \times n$ but 
$2n \times n$), convergence is faster.
\begin{figure}[h!]
\centering
\includegraphics[scale=0.36]{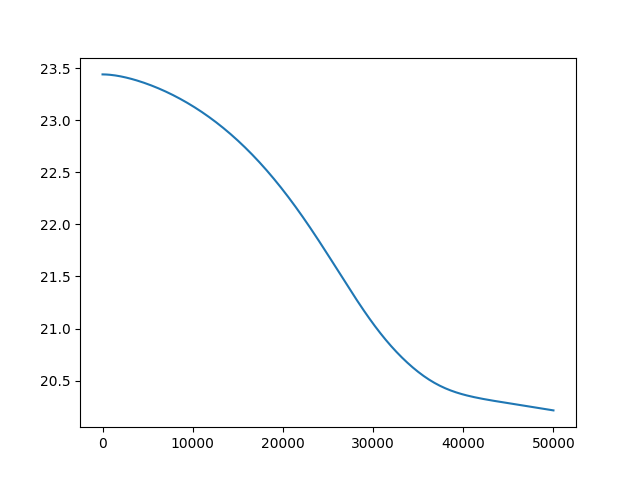}
\includegraphics[scale=0.36]{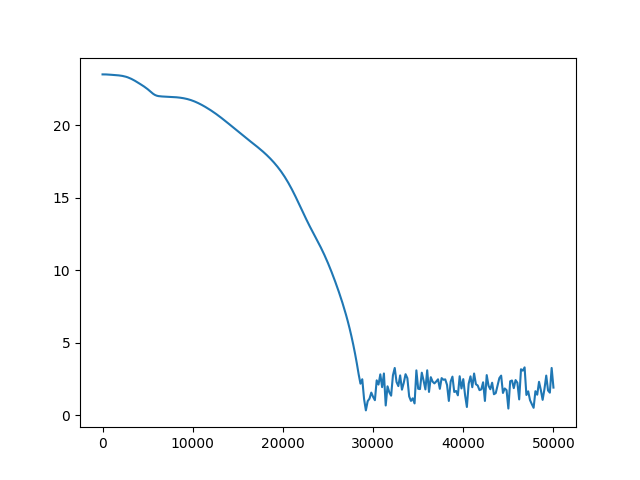}\\
\includegraphics[scale=0.36]{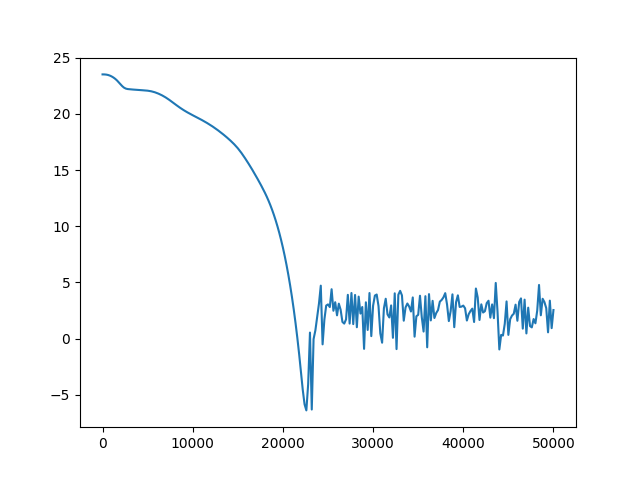}
    \caption{Top left: convergence in the case of many perturbed coefficients, $A=I$. Top right: same, with general $A$. Bottom: same, with a larger $A$ ($2n \times n$).}
\end{figure}
\subsection{"Fine structure": the accuracy in fitting polynomials with many very small coefficients}
In these experiments, we tested the accuracy of fitting a
SOS polynomial with a very simple underlying structure, 
but with all its coefficients perturbed by small random noise. The polynomials were defined by
\[
10\left(\sum_{i=1}^n x_i^2\right)^2 + r(x_1\ldots x_n)
\]
where $r(x_1\ldots x_n)$ is a polynomials all of whose
coefficients are random uniform in $[-0.01,0.01]$. Our
approach yielded a total squared error of $5.5\cdot10^{-10}$, which corresponds to an RMS error of $1.1 \cdot 10^{-7}$ per
coefficients, far smaller than the average
magnitude of the "noisy" coefficients (0.0029).
\begin{figure}[h!]
\centering
\includegraphics[scale=0.36]{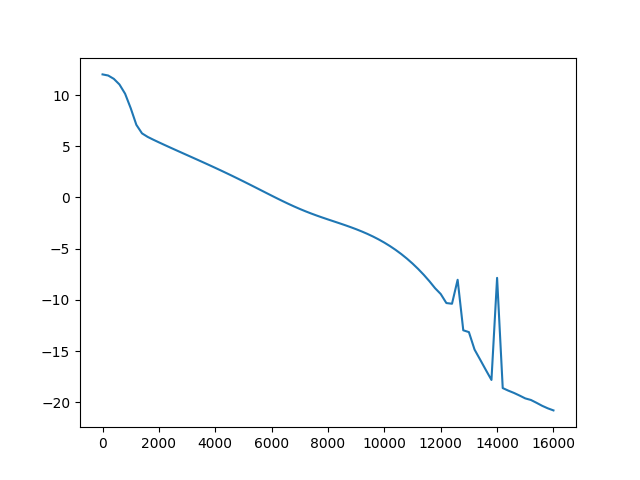}
    \caption{Convergence for a polynomials with many small
    random coefficients.}
\end{figure}
\section{Acknowledgment}
We are grateful to Nicolas Boumal, Richard Zhang, and 
Chenyang Yuan for their helpful advice.

\newpage
\bibliographystyle{acm}
\bibliography{bib}

\end{document}